\documentclass{amsart}
\usepackage{a4wide}
\usepackage[latin1]{inputenc}
\usepackage[T1]{fontenc}
\usepackage{amsmath,amssymb,amsthm,stmaryrd}
\usepackage{mathrsfs}
\usepackage{esint}
\usepackage{color}
\usepackage{pdfsync}

\swapnumbers

\numberwithin{equation}{section}

\makeatletter
  \let\c@subsection\c@equation
\makeatother

\theoremstyle{plain}

\theoremstyle{definition}

\theoremstyle{remark}

\makeatletter
\@namedef{subjclassname@2010}{
  \textup{2010} Mathematics Subject Classification}
\makeatother

\title[Wavelets in spaces of homogeneous type]{addendum to : Orthonormal bases of regular wavelets in spaces of homogeneous type}
\author{Pascal Auscher}
\address{Laboratoire de Math\'ematiques, UMR 8628, Univ. Paris-Sud and CNRS,   F-91405 {\sc Orsay}} 
\email{pascal.auscher@math.u-psud.fr}
\author{Tuomas Hyt\"onen*}\thanks{*Corresponding author. Telephone: +358 2941 51430.}
\address{Department of Mathematics and Statistics, P.O. Box 68, FI-00014 University of Helsinki, Finland}
\email{tuomas.hytonen@helsinki.fi}


\begin{document}

\begin{abstract}
We bring a precision to our cited work concerning the notion of ``Borel measures'', as the choice among different existing definitions impacts on the validity of the results.

\end{abstract}

\maketitle

We wish to bring a precision to our work \cite{AH}. The same remarks apply to the follow-up article~\cite{HT}.  Proposition 4.5 in \cite{AH} states that if $\mu$ is  a non-trivial Borel measure on a quasi-metric space $X$, finite on bounded Borel sets, and $1\leq p<\infty$, then H\"older-$\eta$-continuous functions of bounded support are dense in $L^p(\mu)$, where $\eta$ is the H\"older exponent of the splines constructed in \cite{AH}.

Depending on the meaning of ``Borel measure'', as different definitions can be found in the literature, the result is correct or wrong.

If $\mu$ is a $\sigma$-additive measure on the Borel $\sigma$-algebra, then the result with the given proof is correct, as the measurable sets coincide with the Borel sets.  
In that case, all our results are valid as stated. 

However, if $\mu$ is an outer measure or a $\sigma$-additive measure on $X$ for which the Borel sets are $\mu$-measurable, then for the first sentence of the proof  to be valid one must add the condition that  for every $\mu$-measurable set $A$ (in the sense of Caratheodory for the outer measure case),  there is a Borel set $B\supseteq A$ such that $\mu(A)=\mu(B)$, and the rest of the proof goes through. In \cite{Ma} for example, Borel outer  measures are called regular if this condition holds for all $A$ (not necessarily $\mu$-measurable).
With such a definition of a Borel measure, this regularity condition should be added to our statements.  Without regularity, the correct conclusion of Proposition 4.5  
is density in the space of $L^p$ functions having  a Borel measurable representative. Thus,  our wavelet representations  are valid for functions in this subspace and $1<p<\infty$. This is enough for many purposes.

\bibliographystyle{acm}

\end{document}